\documentclass[10pt]{article}

\usepackage[margin=1in]{geometry}

\usepackage{url,tabu,booktabs,multirow,amsmath,amssymb,amsthm,graphicx,cite}
\graphicspath{.}

\theoremstyle{plain}

\theoremstyle{definition}

\theoremstyle{remark}
\newtheorem{remark}{Remark}

\usepackage{xcolor}

\def\Dprime{{D_\varepsilon}}
\def\vu{{\mathbf u}}

\begin{document}

\begin{center}
\begin{Large}

{\bf Particle Simulation of Fractional Diffusion Equations} \\
\end{Large}

\bigskip

S. Allouch$^1$, M. Lucchesi$^1$, O.P. Le Ma\^{\i}tre$^2$, K.A. Mustapha$^3$, O.M. Knio$^1$\footnote{omar.knio@kaust.edu.sa} \\
$^1$ King Abdullah University of Science and Technology, Thuwal, Saudi Arabia \\
$^2$ CNRS, LIMSI, Universit\'e de Paris Saclay, Orsay, France \\
$^3$ King Fahd University of Petroleum and Minerals, Dhahran, Saudi Arabia 

\end{center}

\begin{abstract}
    This work explores different particle-based approaches to the simulation of one-dimensional fractional sub-diffusion equations in unbounded domains.  
    We rely on smooth particle approximations, and consider four methods for estimating the fractional diffusion term.  
    The first method is based on a direct differentiation of the particle representation; it follows the Riesz definition of the fractional derivative and results in a non-conservative scheme.  
    The other three methods follow the particle strength exchange (PSE) methodology and are by construction
    conservative, in the sense that the total particle strength is time invariant.  
    The first PSE algorithm is based on using direct differentiation to estimate the fractional diffusion flux, and exploiting the resulting
    estimates in  an integral representation of the divergence operator. 
    Meanwhile, the second one relies on the regularized Riesz representation of the fractional diffusion term to derive a suitable interaction formula acting directly on the particle representation of the diffusing field.  
    A third PSE construction is considered that exploits the Green's function of the fractional diffusion equation. 
    The performance of all four approaches is assessed for the case of a one-dimensional diffusion equation with constant diffusivity.  
    This enables us to take advantage of known analytical solutions, and consequently conduct a detailed analysis of the performance of the methods.  
    This includes a quantitative study of the various sources of error, namely filtering, quadrature, domain truncation, and time integration, as well as a space and time self-convergence analysis.  
    These analyses are conducted for different values of the order of the fractional derivatives, and computational experiences are used to gain insight that can be used for generalization of the present constructions.
\end{abstract}


\section{Introduction}

A wide variety of transport problems of relevance to industrial and environmental 
applications involve anomalous diffusion phenomena.  A canonical example concerns
the well-known problem of the space fractional diffusion equation, which in its simplest
form is described in terms of a fixed, spatially-uniform diffusivity coefficient and a fixed 
fractional order derivatives.  Numerous studies have focused on the numerical solution
of this problem, using both finite-difference and variational methodologies; see 
e.g.~\cite{LynchEtAl2003, TadjeranMeerschaertScheffler2006, ZhangLiuAnh2010,  LiXu2010, 
Mustapha2011, CelikDuman2012, Mustapha2013, DengHesthaven2013,  Mustapha2014, TianZhouDeng2015, Mustapha2015, LiuaYanKhan2017, ZhouTianDeng2013} and references therein.
In contrast to grid-based and Galerkin machineries, experiences with particle-based methods 
are scarce, and appear to be limited to random-walk
type approaches; see e.g.~\cite{Gorenflo2002,Gorenflo2004,Gorenflo2007,Zhang2007,Shi2013,Luchko2016}.

The present work aims at exploring the application of so-called smooth particle approximations~\cite{Cottet:2000} 
to space fractional diffusion equations.  It is specifically motivated by the desire to capitalize on some of the
advantages offered by particle schemes, particularly their natural ability to deal with situations involving unbounded
domain, and to treat far-field boundary conditions.  

To our knowledge smooth particle methods have not yet been applied to fractional diffusion
equations.  Consequently, in this work we focus on developing different smooth particle constructions 
and examining their properties. As outlined in section~\ref{sec:formulation}, we restrict our attention to 
one-dimensional sub-diffusion equations in an unbounded domain.  Starting from a general statement of the conservation 
law in divergence form, we provide expressions of the fractional diffusion flux, and further restrict 
ourselves to the case of a constant, spatially-uniform diffusivity, which enables us to exploit available
expressions of the fundamental solution~\cite{mainardi2001fundamental}.  Different particle schemes
are then developed in section~\ref{sec:particle approx}.  Specifically, we introduce a smooth particle 
representation of the diffusing field, and use this representation to derive expressions of the fractional
diffusion flux.  Based on this development, we explore two avenues for construction particle schemes.
In the first avenue, the diffusion operator is directly differentiated resulting in a non-conservative approach 
that we label Direct Differentiation (DD).  To enforce conservation properties, the particle strength exchange 
(PSE) formalism is considered in the second avenue.  Following this formalism, we develop three conservative
formulations.  The first PSE formulation is constructed based on approximating the derivative (divergence) of 
the fractional diffusion flux; this approach is referred to as Flux-PSE (FPSE).  Another alternative is explored 
in which the PSE approximation is applied the second derivative of a convolution integral, leading to what we
refer to as the Kernel PSE (KPSE) method.   Lastly, we propose a third PSE variant based on the Green's function 
of the problem, therefore coining it the GPSE scheme.
In section~\ref{sec:res}, we analyze the performance of the different
constructions and assess their relative merits.  In particular, a computational study is conducted
of the impact of discretization and truncation parameters.  Major conclusions are summarized in 
section~\ref{sec:conc}.

\section{Formulation}\label{sec:formulation}

We focus on one-dimensional (1D) fractional diffusion equations resulting from the 1D continuity equation:
\begin{equation}
\frac{\partial u}{\partial t} = - \frac{\partial Q^\beta}{\partial x}
\label{eq:conservation}
\end{equation}
where $Q^\beta$ is the fractional diffusion flux of order $\beta$.

For $0 < \beta < 1$, the fractional diffusional flux is first defined according to:
\begin{equation}
Q^\beta(x,t) \equiv - {\mathcal D} \frac{1}{2\sin\left(\frac{\beta\pi}{2}\right)} \left[\frac{\partial ^{ \,\beta} }{\partial x^{\beta}}u(x,t)-\frac{\partial ^{ \,\beta} }{\partial (-x)^{\beta}}u(x,t) \right],
\label{eq:Qdef}
\end{equation}
where ${\mathcal D}$ is the (generally variable but in our case constant) diffusivity, and 
\begin{align}
\frac{\partial^{\,\beta}}{\partial   x^{\beta}}u(x,t) &\equiv\frac{ 1}{\Gamma(1-\beta)}\frac{\partial}{\partial x}\int_{-\infty}^x \frac{u(\xi,t)}{(x-\xi)^{\beta}}\, \mathrm{d}\xi, \label{eq:lder}\\
\frac{\partial^{\,\beta}}{\partial(-x)^{\beta}}u(x,t) &\equiv\frac{-1}{\Gamma(1-\beta)}\frac{\partial}{\partial x}\int_x^{+\infty} \frac{u(\xi,t)}{(\xi-x)^{\beta}}\, \mathrm{d}\xi, \label{eq:rdir}
\end{align}
are the left-sided and right-sided Riemann-Liouville fractional derivatives~\cite{podlubny1998fractional}, respectively.
Combining the left-derivative and right-derivative terms, we can rewrite $Q^\beta$ as:
\begin{equation}
    Q^\beta(x,t)=-\mathcal{D} c_\beta \frac{\partial}{\partial x}\int_{-\infty}^{+\infty}\frac{u(\xi,t)}{\left|x-\xi\right|^{\beta}}\, \mathrm{d}\xi ,
\label{definition flux}
\end{equation}
where
\begin{equation}\label{eq:cbeta-def}
c_\beta \equiv \frac{1}{ 2\Gamma(1-\beta) \sin(\beta\pi/2) }
\end{equation}

Without loss of generality we set $\mathcal{D} = 1$ in the following.
In this case, the conservation law~\eqref{eq:conservation} reduces to the fractional diffusion equation~\cite{mainardi2001fundamental}:
\begin{equation}
\frac{\partial u}{\partial t}={_xD^{\alpha}_0} u, \label{eq:fracdiffeq}
\end{equation}
where $\alpha\equiv\beta+1$ and
\begin{equation}
_xD^{\alpha}_0 u(x,t) \equiv \frac{\Gamma(1+\alpha)}{\pi \alpha(\alpha-1)}\sin\left(\frac{\alpha\pi}{2}\right)\frac{\partial ^{ 2} }{\partial x^{2}} \int_{-\infty}^{+\infty}\frac{u(\xi,t)}{\left|x-\xi\right|^{\alpha-1}}\, \mathrm{d}\xi.
\label{def fractional derivative}
\end{equation}
is the Riesz fractional derivative of order $\alpha$, $1 < \alpha < 2$.

We shall focus on the numerical solution of Eq.~\eqref{eq:fracdiffeq} using different 
smoothed particle approximations.  For this setting, closed form expressions for its fundamental solution are available~\cite{mainardi2001fundamental}, 
as summarized in~\ref{app:fundamental}.  In particular, we take advantage of these expressions to assess the performance
of the various approaches considered.

\section{Particle approximation}\label{sec:particle approx}

In the framework of smooth particle methods, the domain is discretized using a finite collection of particles, defined by their positions, $x_i$, and volumes, $V_i$.
The particle representation of a function, $u(x)$, is expressed as~\cite{Chorin:1973,Leonard:1980,Leonard:1985,Cottet:2000}:
\begin{align}
u(x)=\sum_{i\in {\mathcal I} } V_i u_i \eta_\epsilon(x-x_i),
    \label{eq:part}
\end{align}
where $u_i$ is the strength of particle $i$, ${\mathcal I} \subseteq \mathbb{Z}$ is the particles index set, $\eta$ is a radial regularization kernel of unit mass, and $\epsilon$ is the associated smoothing parameter.
The function $\eta_\epsilon$ is defined according to:
\begin{align}
\eta_\epsilon(x)=\frac{1}{\epsilon^d}\eta\left(\frac{x}{\epsilon}\right),
    \label{eq:etaeps}
\end{align}
where $d$ is the number of spatial dimensions (which is one in our case).
Clearly, $\eta_\epsilon$ tends to the Dirac measure in the sense of distributions as $\epsilon \to 0$.
We rely on the second-order, squared exponential kernel in 1D~\cite{beale1985high,eldredge2002general}:
\begin{align}
\eta(x)=\frac{1}{\sqrt{\pi}} \exp\left( -x^2 \right).
    \label{eq:etadef}
\end{align}

\subsection{Direct Differentiation}
\label{sec:direct}

\begin{figure}[htb]
    \centering\includegraphics[width=0.95\textwidth]{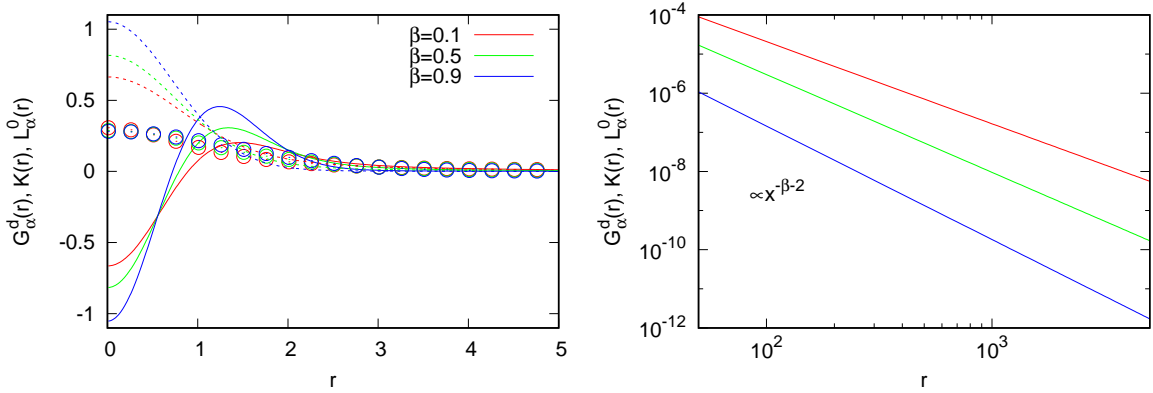}
    \caption{Kernels of the proposed particle schemes: $G^d_\alpha(r)$ from~\eqref{eq:kernel_dd} for the direct differentiation approach (solid lines), $K(r)$ from~\eqref{eq:Ktemp} for the regularized Riesz approach (dashed lines), and $L^0_\alpha(r)$ from~\eqref{eq:l0_forw} and~\eqref{eq:l0_asym} for the Green's function approach (symbols). 
Plotted are curves for three different values of $\beta=0.1$, 0.5, and 0.9, as indicated. Left plot: small values of the argument. Right plot: large values of the argument to illustrate the dependence of the decay rate on $\beta$; only the case of $G^d_\alpha(r)$ is shown for clarity.}\label{fig:kernels_1}
\end{figure}

As implied by its name, in the direct differentiation approach~\cite{Fishelov:1990} one directly differentiates the particle representation of $u$, and then uses a collocation (midpoint rule) approximation to estimate the discrete particle fluxes.
As outlined in~\ref{app:flux}, the particle representation in~\eqref{eq:part} may be readily exploited to estimate the Riemann-Liouville diffusion flux $Q^\beta$, $0 < \beta < 1$.
Specifically, for unit diffusivity we have:
\begin{equation}\label{eq:RLflux}
    Q^\beta(x)=-\frac{2^{\frac{\beta-3}2}}{\sqrt\pi\epsilon^\beta \sin\left(\frac{\beta\pi}2\right)}\sum_{i\in {\mathcal I}}\,V_i u_i \frac\partial{\partial x} \left[ \exp\left(-\frac{X_i^2}4\right)\left[D_{\beta-1}(-X_i)+D_{\beta-1}(X_i)\right] \right] , 
\end{equation}
where  $X_i \equiv \sqrt2(x-x_i)/\epsilon$, and $D_\nu$ denotes the parabolic cylinder function, see~\ref{app:parabolic}.
To evaluate the derivatives, we make use of the relation~\cite{parcylfuncs}:
\begin{align}
    \frac\partial{\partial x} \left[ \exp\left(-\frac{x^2}4\right)D_{\nu}(x) \right]  &=-\exp\left(-\frac{x^2}4\right)D_{\nu+1}(x) ,\label{eq:expdnx2_der}
\end{align}
to obtain:
\begin{equation}
    Q^\beta(x)=-\frac{2^{\frac{\beta-2}2}}{\sqrt\pi\epsilon^{\beta+1}  \sin\left(\frac{\beta\pi}2\right)}\sum_{i\in {\mathcal I}}\,V_i u_i \exp\left(-\frac{X_i^2}4\right)\left[D_{\beta}(-X_i)-D_{\beta}(X_i)\right].
    \label{eq:Qbfield}
\end{equation}
Similarly, we differentiate the above expression to obtain:
\begin{equation}
    \frac{\partial Q^\beta(x)}{\partial x}=-\frac{2^{\frac{\beta-1}2}}{\sqrt\pi\epsilon^{\beta+2} \sin\left(\frac{\beta\pi}2\right)}\sum_{i\in {\mathcal I}}\,V_i u_i \exp\left(-\frac{X_i^2}4\right)\left[D_{\beta+1}(-X_i)+D_{\beta+1}(X_i)\right].
    \label{eq:approx_flux}
\end{equation}
This leads to the following discrete particle approximation of~\eqref{eq:fracdiffeq},
\begin{equation}
    \frac{\partial u_i}{\partial t} =\frac{1}{\epsilon^\alpha} \sum_{j\in {\mathcal I}}\,V_j u_j G^d_\epsilon \left( x_i - x_j \right),
    \label{eq:dd}
\end{equation}
where $G^d_\epsilon$ is the radial kernel given by:
\begin{equation}\label{eq:kernel_dd}
    G^d_\epsilon (r)= \frac{1}{\epsilon}G^d_\alpha\left(\frac{r}\epsilon\right);
    \quad G^d_\alpha(r)=-\frac{2^{\frac{\alpha-2}2}}{\sqrt\pi \cos\left(\frac{\pi\alpha}2\right)} S^{\alpha+1}(r),
\end{equation}
with $S^\nu(r)$ given by~\eqref{eq:Sdef}, see~\ref{app:parabolic}.
The function $G^d_\alpha(r)$ is represented in Fig.~\ref{fig:kernels_1}.  
Note that, similar to the case of a Fickian diffusion, the direct differentiation approach does not lead to a conservative scheme.

\subsection{Particle Strength Exchange}
\label{sec:pse}

This section explores different approaches to constructing a conservative particle scheme based on the PSE 
formalism~\cite{Raviart:1985,Degond:1989,Knio:1992,Cottet:2000,eldredge2002general}.

\subsubsection{Riemann-Liouville Treatment}

\begin{figure}[htb]
    \centering\includegraphics[width=0.95\textwidth]{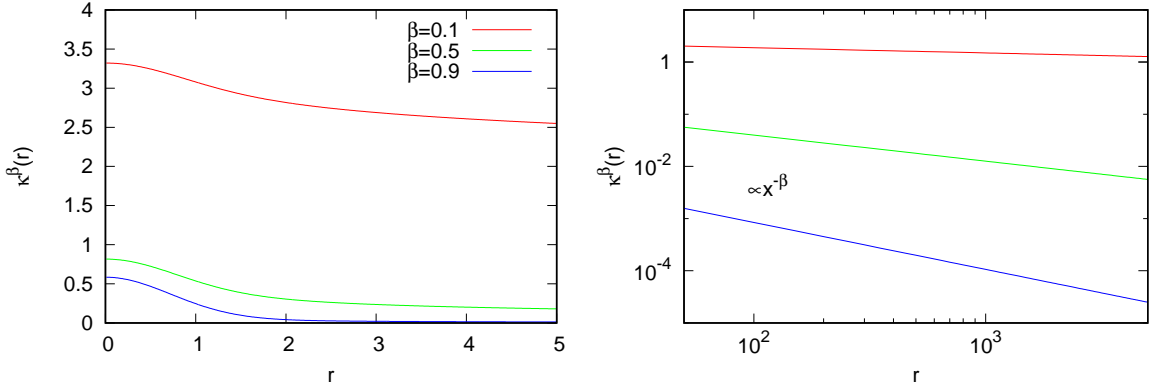}
    \caption{$\kappa^\beta(r)$ (Eq.~\eqref{eq:kernel_utilde}) for three different $\beta=0.1$, 0.5 and 0.9, at small (left) and large (right) values of the argument $r$.\label{fig:kernels_3}}
\end{figure}

According to \eqref{eq:fracdiffeq} and \eqref{def fractional derivative}, for constant diffusivity the fractional diffusional term involves the second derivative of the function
\begin{equation}
\tilde{u}(x,t) \equiv  c_\beta \int_{-\infty}^{+\infty}\frac{u(\xi,t)}{\left|x-\xi\right|^\beta}\, \mathrm{d}\xi,
\label{eq:utilde_def}
\end{equation}
where $\beta = \alpha - 1$ as before.
With the function $u$ expressed as in~\eqref{eq:part}, one can readily evaluate the convolution in~\eqref{eq:utilde_def}; see \ref{app:flux}.
The result can be expressed as:
\begin{equation}
    \tilde{u}(x,t) = \frac{1}{\epsilon^{\beta-1}}  \sum_{i\in {\mathcal I}} V_i u_i \kappa^\beta_\epsilon (x - x_i),
\label{eq:utilde_exp}
\end{equation}
where 
\begin{equation}
\kappa^\beta_\epsilon (r)= \frac{1}{ \epsilon} \kappa^\beta \left( \frac{r}{\epsilon} \right), \quad \mbox{with }
\kappa^\beta(r)= \frac{2^{\frac{\beta-3}{2}}}{\sqrt{\pi}\sin\left(\frac{\beta\pi}2\right)}  S^\beta(r),
\label{eq:kernel_utilde}
\end{equation}
is the radial kernel associated to the second-order kernel $\eta_\epsilon$ in~\eqref{eq:etadef}. 
Figure~\ref{fig:kernels_3} shows $\kappa^\beta$ for different values of $\beta$ and illustrates its slow decay in 
$r^{-\beta}$ when $r\to \infty$, as further discussed below.

Our first attempt at formulating a PSE algorithm consisted in (i)~interpolating $\tilde u$ on the same particle grid used to represent $u$, and (ii)~applying the well-known PSE approximation of the Laplacian~\cite{Raviart:1985,Degond:1989,Knio:1992,Cottet:2000,eldredge2002general} using the exponentially decaying kernel, $\eta$, in~\eqref{eq:etadef}.
This results in the PSE scheme:
\begin{equation}
    \frac{\partial u_i}{\partial t}=\frac{2}{\epsilon^2}\sum_{j\in {\mathcal I}} V_j (\tilde u_j-\tilde u_i)\Phi_\epsilon(x_j-x_i)
    \label{eq:dpse}
\end{equation}
where
\begin{equation}
    \Phi_\epsilon(r)= \frac{1}{\epsilon} \Phi\left(\frac{r}{\epsilon} \right), \quad \mbox{with } \Phi(r)= - \frac{1}{r} \frac{d \eta}{dr} = \frac{2}{\sqrt{\pi}} \exp( -r^2).
    \label{eq:kernel_dpse}
\end{equation}
Note that because the kernel $\Phi_\epsilon$ is radially symmetric the method is conservative in the sense that the total particle strength is time-invariant, i.e.\
\begin{equation}
\frac{d}{dt} \left[ \sum_{i\in {\mathcal I}} u_i V_i \right] = 0.
\end{equation}
However, the method requires two passes over the set of particles, namely for the computation of the function 
$\tilde u$ and then for its differentiation.

Because the characteristic length scale of $\tilde{u}$ is evidently larger than that of $u$, especially for the smallest 
values of $\beta$, large approximation errors occur at the edges of the particle grid.
This issue is due to the slow decay of $\kappa^\beta(x)$; specifically, as can be seen from~\ref{app:parabolic},
\begin{equation}
\kappa^\beta(x) \sim c_\beta  x^{-\beta} \quad \mbox{as $x \to \infty$.}
\label{eq:kappa-asym}
\end{equation}
This problem could be tackled by considering a different particle grid to represent $\tilde{u}$, or by seeking an alternative (potentially stretched) kernel that would be better adapted to the distribution of $\tilde{u}$.
Such approaches are outside the scope of the present work, and are left for a future study.

\begin{figure}[htb]
    \centering\includegraphics[width=0.95\textwidth]{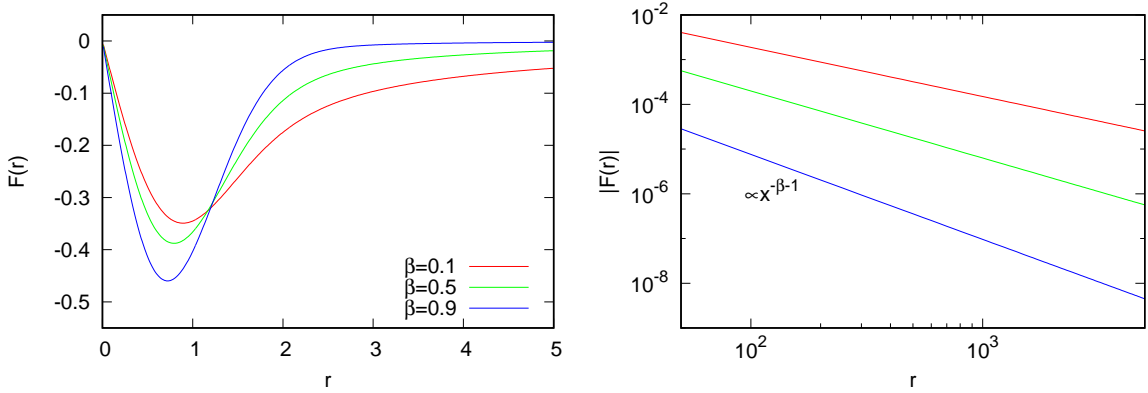}
    \caption{$F(r)$ (see Eq.~\eqref{eq:kernel_fpse}) for different values of $\beta=0.1$, 0.5 and 0.9, and at small (left) and large (right) values of the argument.\label{fig:kernels_2}}
\end{figure}

A possible means of mitigating the issue above is to interpolate $Q^\beta$ on the particle grid, and then apply a symmetrized formula for compute its divergence. 
Using the kernel $\eta$ in~\eqref{eq:etadef} leads to the following expression for the discrete fractional diffusion flux at the particle positions,
\begin{equation}
    Q^\beta_i =-\frac{1}{\epsilon^{\beta}}\sum_{j\in {\mathcal I}}\,V_j u_j F_\epsilon(x_i -x_j),
\end{equation}
where
\begin{equation}\label{eq:kernel_fpse}
    F_\epsilon(r)=\frac1\epsilon F(\frac r\epsilon); \quad F(r)= \frac{2^{\frac{\beta-2}2}}{\sqrt\pi  \sin\left(\frac{\beta\pi}2\right)} 
    T^\alpha(r) ,
\end{equation}
and
\begin{equation}
T^\nu(z) \equiv e^{\frac{-z^2}{2}}\left(D_{\nu-1}(-\sqrt{2}z) - D_{\nu-1}(\sqrt{2}z)\right), \label{eq:Tdef}
\end{equation}
The radial function $F(r)$ is plotted in Fig.~\ref{fig:kernels_2}. It is seen that the function is negative, that it has
a single minimum at $r=0$, and that it increases slowly to $0$ as $r\to \infty$.

We may now estimate the divergence of $Q^\beta$ using~\cite{eldredge2002general}:
\begin{equation}
\left. \frac{\partial Q^\beta}{\partial x} \right|_{x_i} = \frac{1}{\epsilon} \sum_{j\in {\mathcal I}} V_j (Q^\beta_j + Q^\beta_i) \eta^1_\epsilon (x_i - x_j).
\end{equation}
Substituting into the governing equations results in the following particle scheme:
\begin{equation}
    \frac{\partial u_i}{\partial t}= - \frac{1}{\epsilon} \sum_{j\in {\mathcal I}} V_j (Q^\beta_j + Q^\beta_i) \eta^1_\epsilon (x_i - x_j).
    \label{eq:fpse}
\end{equation}
In the computations, we rely on the first-derivative, full-space, second-order kernel~\cite{eldredge2002general},
\begin{equation}\label{eq:kernel_eta1}
    \eta^1_\epsilon (r) = \frac1\epsilon\eta^1\left(\frac r\epsilon\right);\quad \eta^1 (r) = -\frac{2r}{\sqrt{\pi}} \exp ( - r^2 )
\end{equation}
to estimate the flux divergence in~\eqref{eq:fpse}.
We observe that similar to the previous scheme in~\eqref{eq:dpse}, the particle approach in~\eqref{eq:fpse} requires two passes over the particles set to estimate successively the flux and its divergence.

\subsubsection{Riesz Treatment}\label{ssec:Riesz}

We also explored relying on the Riesz representation of the fractional derivative to construct a PSE scheme that acts directly on the particle 
distribution, $u$.

For $1 < \alpha < 2$, provided the function $f$ has continuous second derivative ($f\in C^2(\mathbb R)$), the Riesz fractional derivative has the regularized representation~\cite{mainardi2001fundamental}:
\begin{equation}
    _xD_0^\alpha f(x) = \frac{\Gamma(1+\alpha)}{\pi} \sin\left(  \frac{\alpha\pi}{2} \right) \int_0^\infty \frac{f(x+\xi) - 2f(x) + f(x-\xi) }{\xi^{1+\alpha}} {\mathrm d}\xi.
 \label{eq:Riesz-0}
\end{equation}
A straightforward change of variable leads to:
\begin{equation}
    _xD_0^\alpha u(x) = \frac{\Gamma(1+\alpha)}{\pi} \sin\left(  \frac{\alpha\pi}{2} \right) \int_{-\infty}^\infty \frac{u(y) - u(x)}{\left| y-x \right| ^{1+\alpha}} {\mathrm d}y.
\label{eq:Riesz}
\end{equation}
This suggests several avenues to constructing PSE schemes, based either on regularizing the kernel in~\eqref{eq:Riesz} or regularizing the diffusing field.

We initially considered injecting a regularized version $u_\epsilon$ of the diffusing field $u$,
\begin{equation}
u(x) \approx u_\epsilon(x) \equiv \int_{-\infty}^\infty u(z) \eta_\epsilon(x-z) dz,
\end{equation}
into \eqref{eq:Riesz-0}.
This leads to the following approximation:
\begin{equation}
    _xD_0^\alpha u(x) \approx {_xD_0^\alpha} u_\epsilon(x) = \frac{1}{\epsilon^\alpha } \int_{-\infty}^\infty u(z) G_\epsilon(x-z) dz,
\end{equation}
where
\begin{equation}    
G_\epsilon(r) \equiv  \frac{\Gamma(1+\alpha)}{\pi} \sin\left(  \frac{\alpha\pi}{2} \right) \epsilon^{\alpha} \int_0^\infty \frac{ \eta_\epsilon(r+\xi) - 2 \eta_\epsilon(r) + \eta_\epsilon(r-\xi) } {\xi^{1+\alpha}} {\mathrm d}\xi .
\label{eq:Gd2}
\end{equation}
Performing the integration, and using the recursion properties of parabolic cylinder functions, we conclude that the kernel in Eq.~\eqref{eq:Gd2} is the same as the kernel $G^d_\alpha$ in Eq.~\eqref{eq:kernel_dd}.  This is not surprising because the kernel $\eta$ is smooth and rapidly decaying,
so that the regularized Riesz and Riemann-Liouville representations coincide.
In other words, one recovers the direct differentiation scheme in Eq.~\eqref{eq:dd}.

If on the other hand a shifted approximation is used, one obtains the ``symmetrized'' approximation:
\begin{equation}
    _xD_0^\alpha u(x) \approx  \frac{1}{\epsilon^\alpha } \int_{-\infty}^\infty (u(z) - u(x)) G^d_\epsilon(x-z) dz,
\end{equation}
which leads to the following PSE scheme:
\begin{equation}
    \frac{\partial u_i}{\partial t} = \frac{1}{\epsilon^\alpha} \sum_{i\in {\mathcal I}} V_j \left( u_j - u_i \right) G^d_\epsilon(x_j - x_i).
\label{eq:pseGd}
\end{equation}
This scheme has the unappealing property that it relies on a  non-monotonic kernel that changes sign and assumes negative values at the origin (not shown).
Computational tests using the scheme in Eq.~\eqref{eq:pseGd} exhibited large errors and signs of instability (not shown).
Consequently, we decided not to pursue this avenue.

As an alternative to the approach above, we sought to regularize the singular kernel in Eq.~\eqref{eq:Riesz-0}, namely by considering a convolution of the form:
\begin{equation}
{\cal I}(\epsilon) \equiv \frac{c}{\epsilon^\alpha} \int_{-\infty}^\infty \left( f(y)-f(x) \right) K_\epsilon(x-y) dy
\end{equation}
where $c$ is a constant to be determine and $K$ is a smooth, radial kernel.
Clearly, we wish to select $K$ and $c$ in such a way that ${\cal I}(\epsilon) \to {_x}D^{\alpha}_0 f$ as $\epsilon \to 0$.
A suitable template for $K$ is:
\begin{equation}
K(r)= -\frac 1r\frac{\partial \kappa^\beta}{\partial r} 
\label{eq:Ktemp}
\end{equation}
where $\kappa^\beta$ is defined as in~\eqref{eq:kernel_utilde}.  A straightforward calculation enables us to 
relate the kernels $K$ and $F$, namely through $K(r) = - F(r)/r$.

To justify this choice, and also determine the constant $c$, we substitute the template~\eqref{eq:Ktemp} into the definition of ${\cal I}(\epsilon)$.
With $\epsilon = 1$, this yields:
\begin{eqnarray}
{\cal I}(1)  =  - c \int_0^\infty \left[ f(x+\xi) - 2 f(x) + f(x-\xi)\right] \frac{1}{\xi} \frac{d \kappa^\beta}{d\xi} d\xi.
\end{eqnarray}
Assuming that $f\in C^2(\mathbb R)$, by integrating by parts we obtain:
\begin{equation}
{\cal I}(1) = -c {\cal M}(1) + c {\cal N}(1),
\end{equation}
where 
\begin{equation}
{\cal M}(1) = \int_0^\infty \frac{f(x+\xi) - 2 f(x) + f(x-\xi)}{\xi^{\alpha+1}} \xi^\beta k^\beta(\xi) d\xi  
\end{equation}
and 
\begin{equation}
{\cal N}(1) = \int_0^\infty \frac{f'(x+\xi) - f'(x-\xi)}{\xi^\alpha} \xi^\beta k^\beta(\xi) d\xi.
\end{equation}
If the argument of $K$ is stretched via a scale factor $\epsilon$, then in the limit $\epsilon \to 0$ the quantity $\xi^\beta k^\beta(\xi)$ will approach a constant which can be readily obtained from~\eqref{eq:kappa-asym}.
Consequently, as $\epsilon \to 0$
\begin{equation}
{\cal M}(\epsilon) \to c_\beta  \left[ \frac{\Gamma(1+\alpha)}{\pi} \sin\left(  \frac{\alpha\pi}{2} \right) \right]^{-1} {_xD^{\alpha}_0 f},
\end{equation}
whereas
\begin{equation}
{\cal N}(\epsilon) \to \alpha c_\beta \left[ \frac{\Gamma(1+\alpha)}{\pi} \sin\left(  \frac{\alpha\pi}{2} \right) \right]^{-1} {_xD^{\alpha}_0 f}.
\end{equation}
Combining the two terms, we have: 
\begin{eqnarray}
{\cal I}(\epsilon) & \to & \frac{c (\alpha-1)\pi}{2 \displaystyle \sin\left(\frac{\pi \beta}{2}\right) \sin\left(  \frac{\alpha\pi}{2} \right) \Gamma(1-\beta) \Gamma(1+\alpha) } {_xD^{\alpha}_0 f} \nonumber \\
 & = & \frac{c (\alpha-1)\pi}{ \sin\left( \pi \beta \right) \Gamma(1-\beta) \alpha (\alpha -1) \Gamma(\beta) } {_xD^{\alpha}_0 f} \nonumber \\
 & = & \frac{c}{\alpha} {_xD^{\alpha}_0 f}.
\end{eqnarray}
The desired behavior (limit) is consequently obtained with $c = \alpha$.

Using the template of $K$ in~\eqref{eq:Ktemp} with $c=\alpha$ results in the following PSE scheme:
\begin{equation}
    \frac{\partial u_i}{\partial t} = \frac\alpha{\epsilon^\alpha} \sum_{j\in {\mathcal I}} V_j \left( u_j - u_i \right) K_\epsilon(x_j - x_i),
\label{eq:kpse}
\end{equation}
where
\begin{equation}\label{eq:kernel_kpse}
    K_\epsilon (r) = \frac{1}{\epsilon} K\left(\frac r\epsilon\right) .
\end{equation}
The function $K(r)$, defined in \eqref{eq:Ktemp}, is plotted in Fig.~\ref{fig:kernels_1}.
Note that the kernel $K$ is positive and exhibits the expected asymptotic behavior as $r \to \infty$.

\begin{remark}
    It is interesting to contrast the PSE constructions for the Fickian ($\alpha  = 2$) and fractional ($1 < \alpha < 2$) diffusion.  
    Starting from the radial, unit mass, rapidly-decaying kernel $\eta$ in~\eqref{eq:etadef}, in the Fickian case we rely on the template $-r^{-1} d\eta / dr$, and convolve the Taylor expansion of $u$ to extract a PSE approximation of the Laplacian.
    In the fractional diffusion case, we start by convolving $\eta$ with the Riemann-Liouville kernel $|r|^{-\beta}$ to first construct the stretched kernel $\kappa^\beta$.
    We then rely on the template $-r^{-1} d\kappa^\beta / dr$ to appropriately remove the singularity of the Riesz representation of the fractional derivative.
    The PSE approximations for both the Fickian and non-Fickian diffusion cases can be captured in a single, order-dependent, expression of the form:
    \begin{equation}
        \frac{\alpha}{\epsilon^\alpha} \sum_{j\in {\mathcal I}} V_j \left( u_j - u_i \right) {\cal L}_\epsilon(x_j - x_i) \quad 1 < \alpha \leq 2,
    \end{equation}
where ${\cal L}_\epsilon$ coincides with $\Phi_\epsilon$ from~\eqref{eq:kernel_dpse} in the Fickian case, and with $K_\epsilon$ from~\eqref{eq:kernel_kpse} for $1 < \alpha < 2$.
Note that for the Fickian case, the diffusion kernel ($\Phi$) has finite second moment; in fact for the present choice of $\eta$ all the moments of $\Phi$ are finite.  On the other hand, for the fractional diffusion case, the kernel ($K$) has finite mass and first moment, but infinite higher moments.
This highlights the singular nature of the limit $\alpha \to 2$.
\end{remark}

\begin{remark}
    We also explored an ad-hoc variant of the approach just presented, based on regularizing the Riesz kernel $\xi^{-\alpha - 1}$ directly using a cutoff lengthscale $\epsilon$.
    In preliminary tests, this classical~\cite{Leonard:1980,Cottet:2000} regularization approach resulted in poor approximations of the diffusion flux, especially at larger values of $\alpha$.
    For brevity, these attempts are omitted from the present discussion.
\end{remark}

\subsubsection{Green's function approach}

The availability of an analytical Green's function provides an alternative means of constructing a conservative strength exchange methodology.
In the constant diffusivity case, the Green's function is given by~\cite{mainardi2001fundamental}:
\begin{equation}
{\cal G}^0_\alpha(x,t) = \frac{1}{t^\gamma} L^0_\alpha \left(   \frac{x}{t^\gamma}  \right),
\end{equation}
where $\gamma \equiv 1/ \alpha$. 
For $\Delta t>0$ and using this Green's function, one can express the change in the solution at an arbitrary point $x$ as~\cite{Choquin:1989}:
\begin{eqnarray}
u(x,t+\Delta t) - u(x,t) & = & {\cal G}^0_\alpha(x,\Delta t) * u(x,t) - u(x,t)  \nonumber \\
                         & = & \int_{-\infty}^\infty {\cal G}^0_\alpha(x-y ,\Delta t)  u(y,t) {\mathrm d}y - u(x,t) \nonumber \\
                         & = & \int_{-\infty}^\infty {\cal G}^0_\alpha(x-y ,\Delta t)  \left[ u(y,t) - u(x,t) \right]{\mathrm d}y.
\end{eqnarray}
Using superscripts $n$ and $n+1$ to refer to the time levels, so that~$\Delta t=t^{n+1}-t^n$, the previous representation readily leads to the following particle approximation:
\begin{equation}
u_i^{n+1} - u_i^n = \sum_{j\in {\mathcal I}} V_j \left( u_j^n - u_i^n \right) E_{\epsilon}(x_j - x_i),
\label{eq:Green-dt1}
\end{equation}
where $E_{\epsilon}$ is the radial kernel given by:
\begin{equation}\label{eq:kernel_gpse}
    E_{\epsilon} (r) = \frac{1}{\epsilon} E\left(\frac r\epsilon\right) \quad\mbox{with } E(r)= L^0_\alpha(r), \mbox{ and }
\end{equation}
and 
\begin{equation} \label{eq:defepsilonGPSE}
    \epsilon = (\Delta t)^\gamma.
\end{equation}
The kernel $L^0_\alpha(r)$ is shown in Fig.~\ref{fig:kernels_1} with symbols. It is seen to be positive with monotonic decay to $0$. 
The particle scheme in~\eqref{eq:Green-dt1} is clearly conservative, and it can be alternatively expressed as:
\begin{equation}\label{eq:gpse}
\frac{u_i^{n+1} - u_i^n}{\Delta t} = \frac{1}{\epsilon^\alpha} \sum_{j\in {\mathcal I}} V_j \left( u_j^n - u_i^n \right) E_{\epsilon}(x_j - x_i),
\end{equation}
which enables us to contrast the present scheme with the PSE formulations above.

\section{Results}\label{sec:res}

\subsection{Reference problem}\label{sec:referencePB}
In this section, the behavior of the particle schemes is investigated, focusing on the effect of the main numerical parameters: the kernels width $\epsilon$, the number of particles $N$ (that is ${\rm Card}(\cal I))$, the half-width $D$ of the discretized domain, and the time step $\Delta t$.
As a consequence of both the slow decay of some of the kernels, and of the characteristic length scale associated to the fundamental solution (see~\ref{app:fundamental}), the particles are initialized in a truncated
domain that is centered at the origin.  The half-width of the domain is denoted by $D$, and is specified in 
terms $R_\alpha$ and the simulation time, $t_f$.  Unless  otherwise noted, we shall use
\begin{equation}
    D=160t_f^{1/\alpha}R_\alpha.
\end{equation}
Thus, the extent of the computational domain varies with the fractional diffusion order, reflecting the 
dependence of the characteristic length on $\alpha$.
We use an odd number of particles uniformly distributed over $[-D,D]$,
with the first and last particle centers lying on the extremes of the domain, whereas one of the particles is centered at $x=0$ to better capture the peak value of the solution.
Although computational efficiency is not the concern of the present work, all the simulations were performed on a 56-core workstation, running a parallel implementation of the code using the shared memory paradigm implemented by the OpenMP library for FORTRAN language.  

The methods considered in the numerical tests are summarized in Table~\ref{tab:methods}.
\begin{table*}[h]
    \centering
    \begin{tabu} to 0.9\textwidth {lX[c]cc}
    \toprule
                                       & referred    & method                                                   & kernel                                             \\
    \toprule
        Direct differentiation         & DD          & Eq.~\eqref{eq:dd}                                        & Eq.~\eqref{eq:kernel_dd}                           \\
    \midrule                                                                                                    
        Flux PSE description           & FPSE        & Eq.~\eqref{eq:fpse}                                      & Eq.~\eqref{eq:kernel_fpse}, \eqref{eq:kernel_eta1} \\
    \midrule                                                              
        Regularized Riesz formulation  & KPSE        & Eq.~\eqref{eq:kpse}                                      & Eq.~\eqref{eq:kernel_kpse}                         \\
    \midrule                                                                                                    
        Green's function based         & GPSE        & Eq.~\eqref{eq:gpse}                                      & Eq.~\eqref{eq:kernel_gpse}                         \\        
    \bottomrule
    \end{tabu}
\caption{Summary of the methods.\label{tab:methods}}
\end{table*}

As reference case, we consider the evolution of the fundamental solution with $\beta=0.5$, from time 
$t_0=0.5$ to $t_f=1.5$.   Following the rule given above, for $\beta=0.5$ we obtain $D=357.5$, and 
we use as a reference discretization a grid with $N=32001$ particles, corresponding to a constant particle volume  $V_j= 2D / N = h \approx 2.23\cdot10^{-2}$. 
For all the methods but GPSE, we select a kernel size $\epsilon=2h$, that is, an overlap ratio of 2, and
we consider a classical first-order Runge-Kutta scheme (RK1) with time step $\Delta t=5\cdot10^{-5}$ to integrate the evolution of the particle strengths.
For the GPSE method, the time step is set to a larger value, $\Delta t=1\cdot10^{-2}$, to obtain a smoothing
parameter $\epsilon$ in Eq.~\eqref{eq:kernel_gpse} which is roughly equal to $2h$, as in the other methods.

\begin{figure}[htb]
    \centering\includegraphics[width=0.95\textwidth]{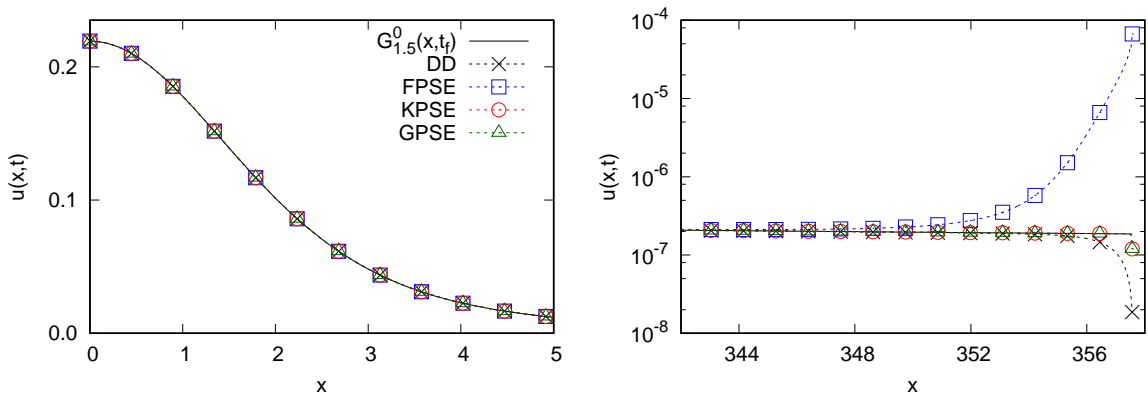}
    \caption{Fundamental solution at $t=1.5$ for $\beta=0.5$ (solid line) compared with the four proposed methods (symbols as indicated). 
    Details on the discretization parameters are provided in the text. Left: solutions for $x\in[0,5]$. Right: zoom on the solutions in the neighborhood of the right boundary of the computational domain.\label{fig:b05_u_ref}}
\end{figure}

The exact solution at $t_f=1.5$, and the four particle solutions are provided in Fig.~\ref{fig:b05_u_ref} for $x\in[0,5]$ (left plot) and $x\in[342,358]$ (right plot). 
It is found that all four methods accurately capture the peak value of the solution and the decay as $|x|$ increases. 
In fact the error is very small for all the methods, except near the boundaries of the computational domain, especially for the DD and FPSE solutions.
It is remarked that the absolute point-wise error remains small everywhere, whereas the relative point-wise 
error increases as $|x|\to D$.
This behavior does not signal instability, but indicates that domain truncation errors are more pronounced
at the computational boundaries.

\subsection{Effect of domain truncation}\label{sec:domconv}
To gain a better understanding of the effect of domain truncation, and how errors incurred at the computational boundary can affect the solution inside the entire domain, we proceed as follows. 
First we focus our attention on the error measured over a subinterval of half-width $\Dprime < D$, and consider the relative $L_1$-error:
\begin{equation}\label{eq:errorL1}
    \varepsilon_{\Dprime}=\frac{\displaystyle\int_{-\Dprime}^{\Dprime} |u(x,t_f)-G_\alpha^0(x,t_f)|\,{\mathrm d}x}{\displaystyle\int_{-\Dprime}^{\Dprime} |G_\alpha^0(x,t_f)|\,{\mathrm d}x}.
\end{equation}
We recall that $G_\alpha^0(x,t_f)$ is the fundamental solution of the reference problem given in~\ref{app:fundamental}.
We set $\Dprime=5R_\alpha$, and we consider domain half widths $D=C t_f^{1/\alpha}R_\alpha$ where the constant $C$ is chosen such that $\Dprime<D$. In the following we use $t_f=1.5$ and $C\in[10,160]$. 
Note that for this definition of $\Dprime$, the solution at $t_f=1.5$, as considered in the following, has roughly 97\% of its mass in the interval $[-\Dprime,\Dprime]$. 
Also note that for $C=160$ one recovers the computational domain of the reference case for $\alpha=1.5$ and the largest value of $D$. 
In order to isolate the effect of the domain truncation, we use a consistent spatial discretization, with the particle volume $h$ depending solely on $\beta$. Specifically, we set $h=2.41\cdot10^{-2}$, $2.23\cdot10^{-2}$, and $1.47\cdot10^{-2}$, for $\beta=0.1$, 0.5, and 0.9 respectively.
The smoothing parameter is set as before to be $\epsilon=2h$; the time-step of the DD, FPSE and KPSE methods is fixed to $\Delta t=5\cdot10^{-5}$ in the RK1 scheme, while the GPSE method uses~\eqref{eq:defepsilonGPSE} and $\epsilon=2h$ to obtain $\Delta t=5\cdot10^{-2},4\cdot10^{-2},1\cdot10^{-3}$ for $\beta=0.1$ ,0.5, and 0.9, respectively.
Finally, the relative error $\epsilon_{\Dprime}$ is approximated by summation over the particles, namely according to:
\begin{equation}
    \varepsilon_{\Dprime} \approx \frac{\displaystyle \sum_{i\in\cal I'} V_i|u_i(t_f)-G_\alpha^0(x_i,t_f)|\,{\mathrm d}x}{\displaystyle\int_{\Dprime}^{\Dprime} |G_\alpha^0(x,t_f)|\,{\mathrm d}x},
\end{equation}
where $\cal I' \subset \cal I$ is the subset of particle indexes such that $|x_i|\le \Dprime$.

\begin{figure}[htb]
    \centering\includegraphics[width=0.95\textwidth]{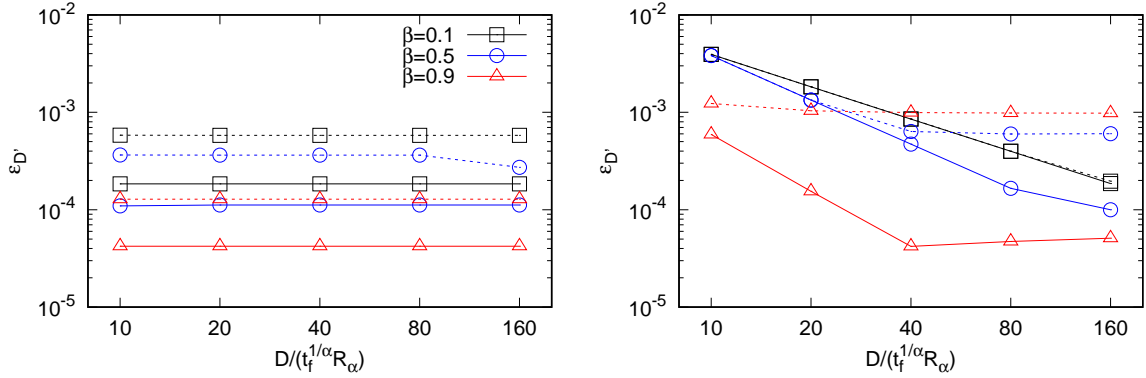}
    \caption{Relative error $\epsilon_{\Dprime}$ versus $D$.  The relative error is defined by~\eqref{eq:errorL1}, and is estimated using 
    $\Dprime= 5R_\alpha$ and the computed solution at $t_f=1.5$.  Curves are plotted for the DD (left, solid line), FPSE (left, dashed line), KPSE (right, solid line) and GPSE (right, dashed line) methods, and three values of $\beta$ as indicated. Other numerical parameters are provided in the 
    text, see section~\ref{sec:domconv}.
    \label{fig:d_relerr}}
\end{figure}

In Fig.~\ref{fig:d_relerr},  the relative error is plotted against $D$, for all four methods considered and for different orders of the fractional derivative.
For the  DD and FPSE methods (left plot), it is seen that with $C=Dt_f^{-1/\alpha}/R_\alpha>10$, the error $\varepsilon_\Dprime$ is essentially independent of $D$ for all values of $\beta$ presented.  This indicates that for $ | x | \le \Dprime$ the error is dominated by the spatial and 
time discretization.  One can also observe that the relative error tends to increase as $\beta$ decreases, and based on this metric the DD method seems more accurate than the FPSE method at the same discretization parameters and value of $\beta$.

Focusing now on the case of the KPSE and GPSE methods, shown in the right plot of Fig.~\ref{fig:d_relerr}, different behaviors are reported as $D$ varies.  First, for the KPSE method (solid lines), an initial decay of $\varepsilon_\Dprime$ with $D$ is observed following which the relative error 
levels off.   The saturation of the relative error occurs at a value of $Dt_f^{-1/\alpha}/R_\alpha$ that increases as $\beta$ decreases. 
In other words, the KPSE method is more sensitive to the truncation of the computational domain when $\beta$ is small. 
This finding is consistent with the analysis of the characteristic length scales of the solution.
Finally, for the GPSE method (dashed lines), one sees that for short computational domains with $C\lesssim 30$
the error is now decreasing as $\beta$ decreases, whereas for sufficiently large $D$ we observe once again that $\varepsilon_\Dprime$ levels off, with higher values for lower $\beta$. This behavior suggests a different mechanism for the effect of the domain truncation in the GPSE method.

In all cases, the results presented justify the value of $D$ that is selected for the reference case, see section~\ref{sec:referencePB}, for which 
domain truncation has negligible impact on $\varepsilon_\Dprime$.  
Finally, it is interesting to remark that, before reaching their respective plateau value as $D$ increases, 
the relative errors $\varepsilon_\Dprime$ of the KPSE and GPSE method seems to follow a decay in $\varepsilon_\Dprime \sim C^{-\alpha}$.
This points to the origin of the domain truncation error, which is related to the decay of the kernels  involved and not due to space or time
discretization parameters.  

\subsection{Effect of space and time discretization parameters}\label{sec:spacetime}

In this subsection, we investigate the  impact of the spatial ($h$) and time ($\Delta t$) discretization parameters of the error incurred by the DD, FPSE and KPSE methods.  
We do not consider here the Green's function based approach (GPSE) because in order to maintain a constant overlap ratio $\epsilon/h$, the values of $h$ and $\Delta t$ can not be varied independently. In the
present numerical experiments, the domain half-width $D$ is kept fixed and coincides with the reference value for $\beta=0.5$, that is $D=357.5$.  We also set $\Dprime = 5R_\alpha$. 

\begin{figure}[htb]
    \centering\includegraphics[width=0.95\textwidth]{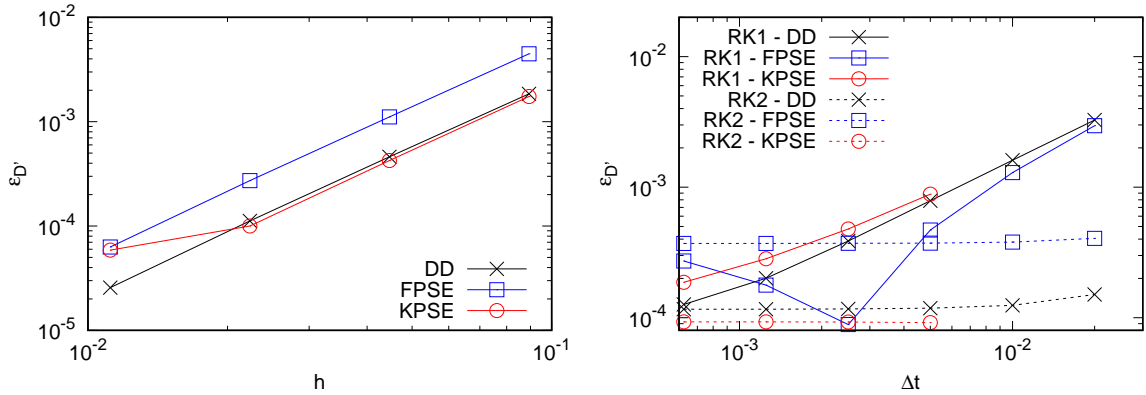}
    \caption{Relative error $\varepsilon_{\Dprime}$ for the DD, FPSE and KPSE methods as a function the particle spacing, $h$ (left), and time step, $\Delta t$ (right).
    The relative error is defined in~\eqref{eq:errorL1}, and is evaluated at $t_f=1.5$ with $\beta=0.5$.
    Other numerical parameters are given in the text, see section~\ref{sec:spacetime}.
    \label{fig:b05_xt_relerr}}
\end{figure}

The left plot in Figure~\ref{fig:b05_xt_relerr} reports $\varepsilon_\Dprime$ for the DD, FPSE and KPSE methods, and for different values of the particles volume $h$. 
In this plot, $h$ is varied by changing the number, $N$, of particles in the computational domain, while the smoothing parameter is adjusted to have $\epsilon=2h$; as before, the RK-1 time-scheme is used with a fixed time step $\Delta t=5\cdot10^{-5}$. 
For the coarsest spatial discretizations (the largest values of $h$), we observe a relative error decreasing as $O(h^2)$ for the three methods, as expected from the moment properties of the kernels. The FPSE method is seen to have the highest spatial error compared to the others at the same $h$.
For the smallest $h$ tested, the error of the KPSE method is seen to level-off as the time discretization error then becomes dominant, because of the fixed time-step. The second-order convergence rate of the relative error $\varepsilon_{\Dprime}$ with $h$ was further verified through additional numerical tests (not shown), which showed the
same trends for $\beta=0.1$ and $\beta = 0.9$.

Next, the effect of the time step is investigated for the reference problem, namely by computing $\varepsilon_{\Dprime}$ for different $\Delta t$ using both RK-1 and RK-2 schemes. The results of these experiments are presented in the right plot of Figure~\ref{fig:b05_xt_relerr}. 
In all these computations, the spatial discretization was the same as in the reference case ($D=357.5$, $N=32,001$, $\epsilon=2h$). 
The curves shows that for the second order Runge-Kutta scheme (RK2, dotted lines), no significant improvement of $\varepsilon_{\Dprime}$ is achieved when $\Delta t$ is lowered, because the time-discretization error is quickly dominated by the spatial discretization error. 
However,  the initial decay of the error for the RK1 scheme is noticeable at the highest range of the time-step values tested, although it subsequently levels off when $\Delta t$ becomes small enough.

To have a better characterization of the spatial and temporal convergence rates of the proposed DD, FPSE and KPSE methods, we estimate 
their orders of self-convergence. Specifically, the self convergence order $p^{(\ell)}$ is estimated from 
\begin{equation}
    p^{(\ell+1)} = \log_2 \displaystyle\frac{\sum_{\in {\cal I}^{(\ell)}}|u^{(l)}_i-u^{(l+1)}_i|}{\sum_{i\in {\cal I}^{(\ell)}}|u^{(\ell+1)}_i-u^{(\ell+2)}_i|}
\end{equation}
where $\ell$ denotes the discretization level and ${\cal I}^{(\ell)}$ the index set of particles belonging to levels $\ell$ to $\ell+2$. 
When going from a level to the next, the numerical parameter investigated is halved. 
Tables~\ref{tab:x_selfconv} and~\ref{tab:t_selfconv} summarizes the order of self-convergence in space and time, respectively, for the 
DD, FPSE and KPSE methods. Regarding the spatial discretization, the expected second-order rate of self convergence is observed. 
For the time-discretization, Table~\ref{tab:t_selfconv} confirms that the convergence in time is governed by the accuracy of the selected 
time-scheme.  Note that these orders of self-convergence are insensitive to $\beta$ and similarly, not affected by (reasonable) changes 
in $\epsilon$.

\begin{table}[hbt]
    \centering
    \begin{tabu} to 0.9\textwidth {cX[c]X[c]X[c]}
        \toprule
        $h$ & DD & FPSE & KPSE \\
        \toprule
        $2.23\cdot10^{-2}$ & 2.00 & 2.01 & 2.00 \\
        \midrule
        $1.11\cdot10^{-2}$ & 2.00 & 2.00 & 2.00 \\
        \bottomrule
    \end{tabu}
    \caption{Summary of the order of self-convergence in space for $\beta=0.5$.\label{tab:x_selfconv}}
\end{table}
\begin{table}[h]
    \centering
    \begin{tabu} to 0.9\textwidth {cX[c]X[c]X[c]X[c]X[c]X[c]}
        \toprule
                   & \multicolumn{2}{c}{DD} & \multicolumn{2}{c}{FPSE} & \multicolumn{2}{c}{KPSE} \\
        \cmidrule{2-7}
        $\Delta t$ & RK1 & RK2 & RK1 & RK2 & RK1 & RK2 \\
        \toprule
        $2\cdot10^{-2}$   & 1.01 & 2.03 & 1.01 & 2.03 & -    & -    \\
        \midrule
        $1\cdot10^{-2}$   & 1.00 & 2.01 & 1.00 & 2.01 & -    & -    \\
        \midrule
        $5\cdot10^{-3}$   & 1.00 & 2.00 & 1.00 & 2.00 & 1.00 & 2.03 \\
        \midrule
        $2.5\cdot10^{-3}$ & 1.00 & 1.98 & 1.00 & 1.98 & 1.00 & 1.92 \\
        \bottomrule
    \end{tabu}
    \caption{Summary of the order of self-convergence in time for $\beta=0.5$.\label{tab:t_selfconv}}
\end{table}

Although not shown, all three methods presently considered exhibit almost no sensitivity to the choice of the smoothing parameter, when the latter
is varied in the range $\epsilon/h\ge2$.

\subsection{Time Stability}

The GPSE method, being based on the Green's function of the problem, is not only conservative but inherently stable, in the sense that an arbitrarily large time-step can be used without compromising the stability of the solution which will remain bounded (in any norm), albeit it may exhibit large error.
However, the other proposed methods, namely DD, FPSE, and KPSE, are generally subject to a stability time-step constraint when an explicit time-scheme is employed. 
It is then interesting to characterize and compare the stability properties of these three methods. 
To this end, we denote $\vu^T = (u_1 \cdot u_N)$ the vector of particles strengths (the superscript $T$ denotes here the transpose).
The evolution of the vector $\vu$ can be written as
\[
    \frac{d \vu}{d t} = A \vu,
\] 
where $A$ is a $N\times N$ matrix whose entries depend on the method considered. Because the particles have fixed positions $x_i$ in the proposed methods, the matrix $A$ is time-invariant. Furthermore, because of the radial nature of the kernels used in these schemes, $A$ is also a symmetric matrix.
Upon introduction of the explicit time discretization, one can derive a \emph{linear} relation between the vectors of particle strengths at time $t_n$ and $t_{n+1}=t_n+\Delta t$, of the form
\[
    \vu^{n+1} = B(\Delta t) \vu^{n},
\]
where the matrix $B(\Delta t)$ is a polynomial in $A$ with coefficients depending on $\Delta t$ and on the selected time discretization.
In the following, we restrict the analysis to the case of the RK1 scheme, which corresponds to 
\begin{equation}\label{eq:matrixtimescheme}
    B(\Delta t) = (I_N +\Delta t A),
\end{equation}
where, $I_N$ is the identity of $\mathbb R^N$. Clearly, $B(\Delta t)$ is symmetric, from the symmetry of $A$, and so has real eigenvalues.
To ensure the stability of the particle scheme, in the sense that $\|\vu^n\|_{L_2} = (\sum_{i} (u^n_i)^2)^{1/2}$ remains bounded for all $n$, 
$B(\Delta t)$ should have all its eigenvalues lying in the unit circle.
Denoting $\lambda_{1\cdots N}$ the eigenvalues of $A$, from~\eqref{eq:matrixtimescheme}  the stability condition becomes
$| 1 + \lambda_i\: \Delta t|\le1$ for $i=1,\cdots,N$. 
Therefore the eigenvalues of $A$ must all verify
\[
    0 \le -\lambda_i \Delta t \le 2.
\]
The condition $\lambda_i\le 0$ is necessary to ensure that the "continuous-in-time" solution $\vu(t) = \exp\left[ A (t-t_0)\right] \vu(t_0)$ does not diverge as $t\to \infty$, and is satisfied for the DD, FPSE and KPSE. 
In fact, the FPSE and KPSE methods being conservative, the largest eigenvalue of the matrix $A$ for these two methods is exactly 0, while one has $\max_{1\le i\le N} \lambda_i <0$ for the DD method which is leaking mass.
Consequently, we only need to focus on the lowest eigenvalue of $A$, denoted $\lambda_{\min}$, and we are led to the stability condition 
$\Delta t \le 2 {|\lambda_{\min}|}^{-1}$.
Defining the nondimensional parameter $a$ as
\begin{equation}\label{eq:stablimit}
    a \equiv\frac{2 {\cal D}}{{|\lambda_{\min}|}h^\alpha},
\end{equation}
the stability condition becomes
\begin{equation}
    \frac{ {\cal D} \Delta t}{h^\alpha}\le a.
\end{equation}
Table~\ref{tab:stablimit} summarizes the value of $a$ for the DD, FPSE and KPSE methods, at different values of $\beta$.
These values were computed for spatial discretization that follows that of the reference case in section~\ref{sec:referencePB} ($N=32,001$ particles), using a power iterations method to compute the leading eigenvalue $\lambda_{\min}$ of $A$. The table shows that the KPSE method has the most stringent stability requirement on the time-step, followed by the DD method, whereas the FPSE method allows a time step that is several times larger
without compromising stability. 
It is also observed that the stability threshold, $a$, does not depend significantly on the value of $\beta$. It is however interesting to note that whereas 
in the DD and FPSE methods the stability limit becomes less stringent as $\beta$ increases, the KPSE method exhibits the opposite trend.
Overall, the stability limit of all the methods are comparable.

\begin{table}[hbt]
    \centering
    \begin{tabu} to 0.9\textwidth {cX[c]X[c]X[c]}
        \toprule
        $\beta$ & DD & FPSE & KPSE \\
        \toprule
        0.1 & 4.81 & 7.05 & 2.25 \\
    \midrule                                                                                                    
        0.5 & 5.25 & 8.83 & 2.17 \\
    \midrule                                                                                                    
        0.9 & 5.43 & 10.5 & 2.04 \\
        \bottomrule
    \end{tabu}
    \caption{Summary of the values of the stability limit $a$, estimated using Eq.~\eqref{eq:stablimit}, for DD, FPSE and KPSE and for different values of $\beta$.\label{tab:stablimit}}
\end{table}


\section{Conclusions}\label{sec:conc}

In this work different particle schemes have been proposed to simulate one-dimensional space-fractional diffusion equations in 
unbounded domains.  This development was motivated by the known reliability and robustness of the particle methods, and their 
ability to accommodate unbounded domains and effectively treat far-field conditions.  Furthermore, we took advantage of the Particle 
Strength Exchange (PSE) formalism possibility to construct conservative particle schemes.

The relative merits of the different schemes can be summarized as follows:
\begin{itemize}
    \item The DD scheme, although nonconservative, has the advantage of being easy to formulate and simple to implement.  In particular, it yields the diffusion flux via single convolution sum over the fields of the individual particles.
    \item On the contrary, the FPSE scheme requires two successive passes over the set of particles to estimate the rate of change of the solution.  However, it is conservative by construction.
    \item Compared to FPSE, the KPSE scheme has the advantage of being computationally less intensive, with the source term obtained in a single convolution sum.  On the downside, the potential extension of the KPSE scheme to the case of variable diffusivity appears to more complicated than in DD and FPSE, which are directly amenable to this generalization.
    \item The GPSE scheme is conservative and has a moderate computational complexity but, as with KPSE, it lacks generality.  In particular, extension of GPSE to the case of variable diffusivity is not straightforward.
\end{itemize} 

Another characteristic of GPSE is that it is based on an exact time integration, owing to the knowledge of the fundamental solution.  This
readily yields  a (fractional order dependent) relation between the time step and the spatial discretization parameter $\epsilon$ (or alternatively the particle volume $h$).  For the other three schemes, namely DD, FPSE, and KPSE, we have performed space and time convergence study to verify that the convergence rate in time indeed corresponds to  the order of the numerical time integration scheme, and that the spatial rate follows the  $O(h^2)$ rate that is expected based on the kernel selected to represent the solution.
We have also shown that this spatial convergence rate is achieved only when the computational domain is sufficiently wide. Otherwise,  the domain truncation error can substantially affect the solution, even at a significant distances from the computational boundaries.  Finally, a brief investigation of the stability of the integration scheme was conducted for the DD, FPSE and KPSE schemes. 
This revealed that all three schemes are subject to comparable time step constraints to ensure stability, which varies weakly with the order of the fractional derivative.

The results of the present work suggest further investigations on the use of particle methods to solve fractional diffusion problems. 
Specifically, one is naturally led to consider a wide class of smooth particle kernels, which would offer the possibility of extending the 
present constructions to achieve higher-order convergence.  An additional avenue that one is encouraged to explore concerns alternative
deterministic particle formulations, for instance based on the diffusion velocity 
approach~\cite{Fronteau:1984,Degond:1990,mycek2016formulation}.  In this context, it appears to be readily possible to capitalize on 
derived expressions for the diffusion flux in building suitable estimates for diffusive particle transport.  The present development also appears
to provide a well-founded foundation for the construction of particle methods for multi-dimensional fraction diffusion problems.  In particular,
multi-dimensional fractional diffusion problems having tensorized structure along spatial directions appear to be addressable by simple extension of the direct differentiation and flux-divergence machineries.

It would also be worthwhile to develop additional capabilities, particularly to enable tackling fractional diffusion problems in complex settings.
An attractive avenue would be the introduction of adaptive particle discretizations, re-meshing strategies, and multiple resolutions approaches. Another avenue concerns the development of fast methods, which offer to dramatically reduce the computational effort needed to estimate
particle interactions involved in flux and source term expressions.  Finally, it would also be interesting to explore stochastic particle 
approaches~\cite{LeMaitre:2007} for problems with uncertain parameters, and the exploitation of related representations to support inference
or calibration problems.

\section*{Acknowledgments}

Research reported in this publication was supported by research funding from King Abdullah University of Science and Technology (KAUST).  

\appendix
\section{Fundamental solution}
\label{app:fundamental}

\setcounter{figure}{0}
\setcounter{table}{0}

\begin{figure}[htb]
    \centering\includegraphics[width=0.65\textwidth]{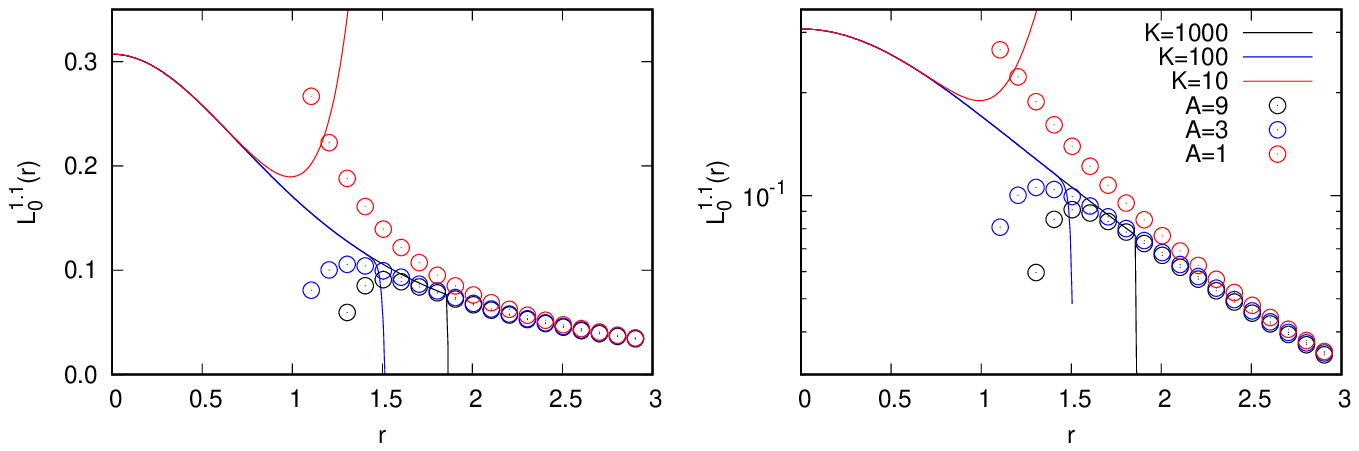}
     \centering\includegraphics[width=0.65\textwidth]{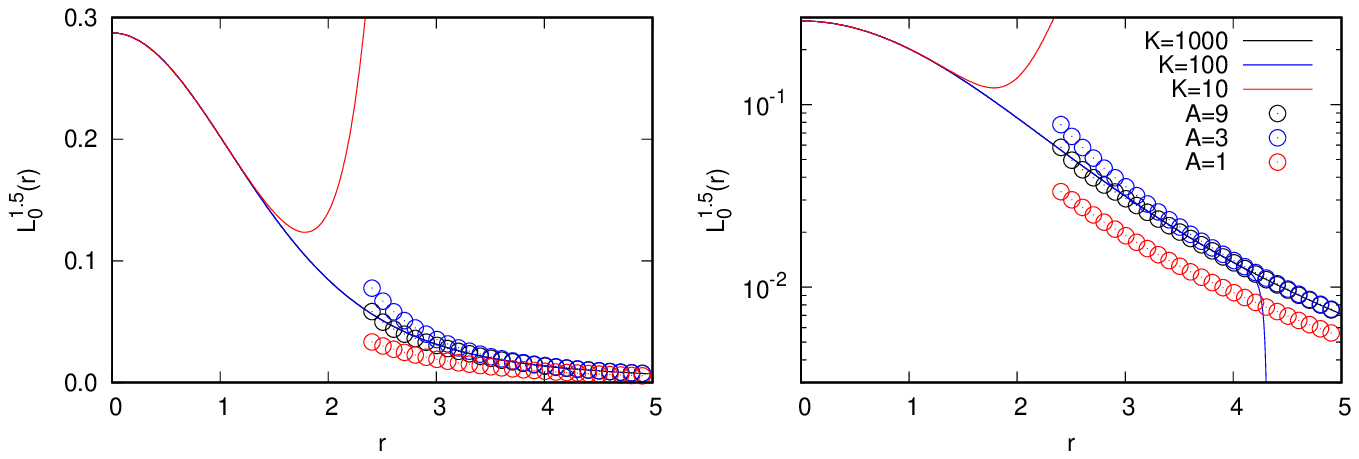}
     \centering\includegraphics[width=0.65\textwidth]{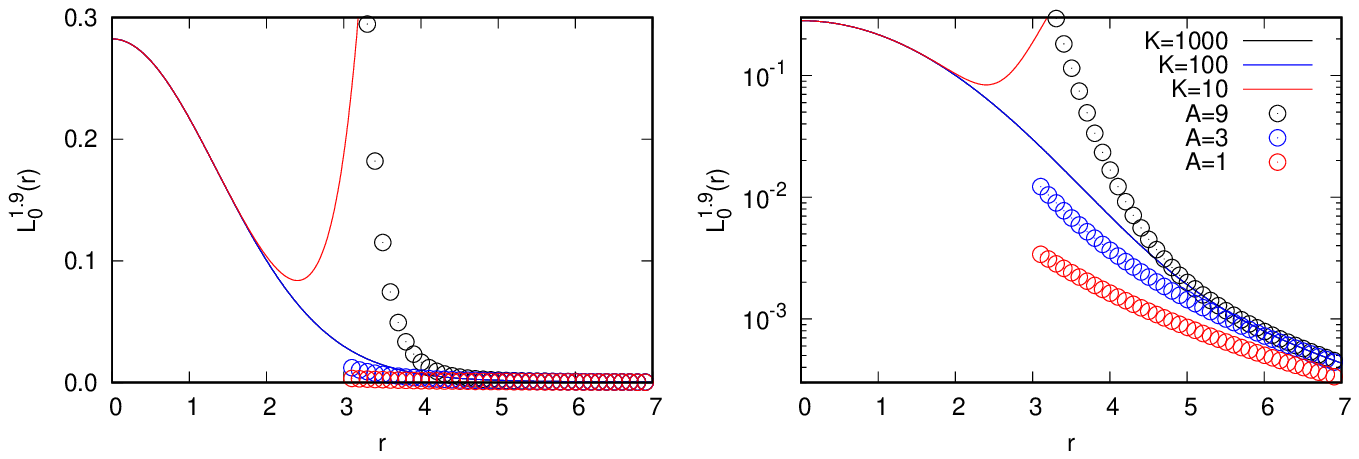}
    \caption{Series expansion of $L^0_\alpha(r)$ (solid lines) and its asymptotic expansion (symbols) represented using respectively $K$ and $A$ terms, as indicated, 
    and for $\alpha=1.1$ (top), $\alpha=1.5$ (center) and $\alpha=1.9$ (bottom).\label{fig:l0_1.1}}
\end{figure}

The Green's function of~\eqref{eq:fracdiffeq} can be expressed as~\cite{mainardi2001fundamental}:
\begin{equation}
G^0_{\alpha}(x,t)=t^{\frac{-1}{\alpha}} L^0_\alpha(x/t^{\frac{1}{\alpha}}),\label{Green function}
\end{equation}
where $L^0_\alpha$ is the so-called reduced Green function.
For $1<\alpha<2$, $L^0_\alpha$ has the series representation:
\begin{equation}
    \label{eq:l0_forw}
L^{0}_\alpha(x)=\frac{1}{x\pi}\sum^{\infty}_{n=1}(-x)^n\frac{\Gamma(1+n/\alpha)}{n!}\sin\left[\frac{-n\pi}{2}\right], \qquad -\infty<x<+\infty.
\end{equation}
Note that $L^0_\alpha$ is radial, i.e.\ $L^{0}_\alpha(-x)=L^{0}_\alpha(x)$,
and that it admits the following asymptotic expansion~\cite{mainardi2001fundamental}:
\begin{equation}
    \label{eq:l0_asym}
L^{0}_\alpha(x)\sim\frac{-1}{x\pi}\sum^{\infty}_{n=1}(-x^{-\alpha})^n\frac{\Gamma(1+ \alpha n)}{n!}\sin\left[\alpha n\frac{\pi}{2}\right], \quad x\to +\infty.
\end{equation}  
Accordingly, for large $|x|$, the Green's function can be approximated by:
\begin{equation}
G_\alpha^0(x,t)\sim\frac{-1}{x\pi}\left[\frac{-t}{x^{\alpha}}\sin(\frac{\alpha\pi}{2})\Gamma(1+\alpha)+\frac{t^2}{x^{2\alpha}}\sin(\alpha\pi)\frac{\Gamma(1+2\alpha)}{2}+\cdots\right], \quad x\to +\infty.
\end{equation}
In Fig.~\ref{fig:l0_1.1}, we use the series representation~\eqref{eq:l0_forw} and asymptotic expansion~\eqref{eq:l0_asym} to illustrate the behavior of $L^{0}_\alpha$ for $\alpha=1.1$, 1.5, and 1.9, respectively.

Because $L^0_\alpha$ does not have a finite second moment, to form a characteristic width of the corresponding distribution we resort to the first moment:
\begin{equation}
    R^\alpha \equiv \frac{ \displaystyle\int_{-\infty}^\infty |x| L^0_\alpha (x) {\mathrm d}x}{ \displaystyle\int_{-\infty}^\infty L^0_\alpha (x) {\mathrm d}x}.
\end{equation}
Noting that the denominator in the above definition is equal to 1, and that $L^0_\alpha (x)$ is radially symmetric, we may rewrite $R^\alpha$ as:
\begin{equation}
    R^\alpha = 2 \int_0^\infty x L^0_\alpha (x) {\mathrm d}x .
\end{equation}
To estimate $R^\alpha$ we decompose the integral according to:
\begin{equation}\label{eq:ra}
R^\alpha = 2 \left(H_1 + H_2 \right),
\end{equation}
where
\begin{equation}
    H_1  \equiv \int_0 ^M  x L^0_\alpha (x) {\mathrm d}x  , \mbox{ and }
    H_2  \equiv \int_M^\infty  x L^0_\alpha (x) {\mathrm d}x \label{eq:H1H2} .
\end{equation}
Substituting~\eqref{eq:l0_forw} and~\eqref{eq:l0_asym} into~\eqref{eq:H1H2}, we obtain:
\begin{equation} \label{eq:h1h2}
    H_1 =\frac M\pi \sum_{k=0}^\infty \frac{(-1)^k\Gamma\left(1+\frac{2k+1}\alpha\right)}{(2k+2)!}M^{2k+1}, \mbox{ and }
    H_2 =\frac M\pi \sum_{n=1}^\infty \frac{(-1)^n\Gamma\left(1+\alpha n \right)}{n!(1-\alpha n)}\sin\left(\frac{\pi \alpha n}2\right)M^{-n\alpha} .
\end{equation}
Estimates of $R^\alpha$ obtained by suitably truncating the above expressions are illustrated in Table~\ref{tab:ralpha}.
\begin{table*}[h]
    \centering
    \begin{tabu} to 0.9\textwidth {cX[c]X[c]X[c]X[c]X[c]X[c]X[c]X[c]X[c]}
    \toprule
    $\alpha$   & 1.1   & 1.2   & 1.3   & 1.4   & 1.5   & 1.6   & 1.7   & 1.8   & 1.9   \\
    \midrule
    $R_\alpha$ & 6.688 & 3.544 & 2.512 & 2.005 & 1.705 & 1.509 & 1.371 & 1.269 & 1.190 \\
    \bottomrule
    \end{tabu}
\caption{Estimates of $R_\alpha$ using Eqs.~\eqref{eq:ra} and \eqref{eq:h1h2}.\label{tab:ralpha}}
\end{table*}

\section{Diffusion flux for $0 < \beta < 1$}
\label{app:flux}

\setcounter{figure}{0}

Substituting~\eqref{eq:part} into~\eqref{definition flux} we get:
\begin{align*}
\frac{-Q^\beta(x)}{c_\beta \mathcal D}=&\frac{\partial}{\partial x}\int_{-\infty}^{+\infty} \frac{\sum_{i=1}^{N} V_iu_i\eta_\epsilon(\xi-x_i)}{|x-\xi|^{\beta}}\, \mathrm{d}\xi\\
                                 =&\frac{\partial}{\partial x}\sum_{i=1}^{N} V_iu_i\int_{-\infty}^{+\infty} \frac{\eta_\epsilon(\xi-x_i)}{|x-\xi|^{\beta}}\, \mathrm{d}\xi\\
                                 =&\sum_{i=1}^{N} V_iu_i\frac{\partial}{\partial x}\int_{-\infty}^{+\infty} \frac{\frac1\epsilon\eta(\frac{\xi-x_i}\epsilon)}{|x-\xi|^{\beta}}\, \mathrm{d}\xi\qquad&\text{see Eq.~\eqref{eq:etaeps}}&\\
                                 =&\sum_{i=1}^{N} V_iu_i\frac{\partial}{\partial x}\int_{-\infty}^{+\infty} \frac1{\epsilon^\beta}\frac{\eta(h)}{|Z_i-h|^{\beta}}\, \mathrm{d}h\qquad &h=\frac{\xi-x_i}\epsilon,\,Z_i=\frac{x-x_i}{\epsilon} . &
\end{align*}
Now letting
\begin{equation}
{\cal K}(x)\equiv \frac{1}{\sqrt{\pi}} \int_{-\infty}^{+\infty} \frac{e^{-h^2}}{|Z_i-h|^{\beta}}\, \mathrm{d}h , \label{eq:calKdef}
\end{equation}
and using Eq.~\eqref{eq:etadef} we have
\begin{align}
\frac{-Q(x)}{\mathcal D}=\frac{c_\beta}{\epsilon^\beta}\frac{\partial}{\partial x}\sum_{i=1}^{N} V_i u_i {\cal K}(x). \label{eq:B1}
\end{align}
Note that the relations above hold for both constant and variable diffusivity, ${\cal D}$.  Next, we decompose the integral in~\eqref{eq:calKdef} according to:
\begin{equation*}\int_{-\infty}^{+\infty} \frac{e^{-h^2}}{|Z_i-h|^{\beta}}\, \mathrm{d}h
=\int_{-\infty}^{Z_i} \frac{e^{-h^2}}{(Z_i-h)^{\beta}}\, \mathrm{d}h +\int_{Z_i}^{+\infty} \frac{e^{-h^2}}{(h-Z_i)^{\beta}}\, \mathrm{d}h.
\end{equation*}
The first integral becomes:
\begin{eqnarray*}
\int_{-\infty}^{Z_i} \frac{e^{-h^2}}{(Z_i-h)^{\beta}}\, \mathrm{d}h
&=& \int_{0}^{+\infty} \frac{e^{-(Z_i-T)^2}}{T^{\beta}}\, \mathrm{d}T\,\,\,\,\,\,\,\,\,\,\,\,\,\,\,\,\,\,\,\,\,\,\,\,(Z_i-h=T)\\
&=&e^{-Z_i^2}\int_{0}^{+\infty}T^{-\beta}e^{-T^2+2TZ_i}\, \mathrm{d}T,
\end{eqnarray*}
whereas the second integral is expressed as:
\begin{eqnarray*}
\int_{Z_i}^{+\infty} \frac{e^{-h^2}}{(h-Z_i)^{\beta}}\, \mathrm{d}h&=& \int_{0}^{+\infty} \frac{e^{-(T+Z_i)^2}}{T^{\beta}}\, \mathrm{d}T\,\,\,\,\,\,\,\,\,\,\,\,\,\,\,\,\,\,\,\,\,\,\,\,(h-Z_i=T)\\
&=&e^{-Z_i^2}\int_{0}^{+\infty}T^{-\beta}e^{-T^2-2TZ_i}\, \mathrm{d}T.
\end{eqnarray*}
Thus,
\begin{eqnarray*}
{\cal K}(x)=\frac{1}{\sqrt{\pi}}e^{-Z_i^2}\Big[\int_{0}^{+\infty}T^{-\beta}e^{-T^2+2TZ_i}\, \mathrm{d}T+\int_{0}^{+\infty}T^{-\beta}e^{-T^2-2TZ_i}\, \mathrm{d}T\Big].
\end{eqnarray*}
For $Re(\vartheta)>0$ and $Re(\beta)>0$, we have~\cite{abramowitz1964handbook}:
\begin{align}
    \label{eq:integral_formula}
\int_{0}^{+\infty}x^{\vartheta-1}e^{-\beta x^2-\mu x}\, \mathrm{d}x=(2\beta)^{\frac{-\vartheta}{2}}\Gamma(\vartheta)e^{(\frac{{\mu}^2}{8\beta})}D_{-\vartheta}(\frac{\mu}{\sqrt{2\beta}}).
\end{align}
For $0<\beta<1$ we use formula~\eqref{eq:integral_formula} to obtain:
\begin{align}
{\cal K}(x)
    =&\frac{e^{-Z_i^2}}{\sqrt{\pi}}\Big[2^{(\frac{\beta-1}{2})}\Gamma(1-\beta)e^{\frac{Z_i^2}{2}}D_{\beta-1}(-\sqrt{2}Z_i)+2^{(\frac{\beta-1}{2})}\Gamma(1-\beta)e^{\frac{Z_i^2}{2}}D_{\beta-1}(\sqrt{2}Z_i)\Big]\nonumber\\
    =&\frac{e^{\frac{-Z_i^2}{2}}2^{(\frac{\beta-1}{2})}\Gamma(1-\beta)}{\sqrt{\pi}}\Big[D_{\beta-1}(-\sqrt{2}Z_i)+D_{\beta-1}(\sqrt{2}Z_i)\Big], \label{eq:frackernel_theorem}
\end{align}
where $Z_i \equiv \frac{x-x_i}{\epsilon}$ and $D_\nu(z)$ is the parabolic cylinder function in Whittaker form.  
Substituting~\eqref{eq:frackernel_theorem} into~\eqref{eq:B1} we get:
\begin{align}
\frac{-Q^\beta(x)}{\mathcal D}=
\frac{2^{(\frac{\beta-3}{2})}}{\sqrt{\pi}\epsilon^\beta \sin(\beta\pi/2) }\frac{\partial}{\partial x}\left[\sum_{i=1}^{N} V_iu_ie^{\frac{-Z_i^2}{2}}\left(D_{\beta-1}(-\sqrt{2}Z_i)+D_{\beta-1}(\sqrt{2}Z_i)\right)\right] \label{eq:genRLflux},
\end{align}
which reduces to~\eqref{eq:RLflux} by setting ${\mathcal D}\doteq 1$.

\section{Parabolic cylinder functions}
\label{app:parabolic}

\setcounter{figure}{0}

The parabolic cylinder function $D_\nu$ is related to the parabolic cylinder functions $U$ and $V$ through the relations~\cite{parcylfuncs}:
\begin{align}
D_{\nu}(z)=U\left(-\frac{1}{2}-\nu,z\right),
\label{eq:DUrel}
\end{align}
and 
\begin{equation}
U(a,-z)=-\sin(\pi a)U(a,z)+\frac{\pi}{\Gamma(\frac{1}{2}+a)}V(a,z).
\label{eq:Uarg}
\end{equation}
The functions $U(a,z)$ and $V(a,z)$ have the following power series representations:
\begin{align}
    U(a,z)&=U(a,0)u_1(a,z)+U'(a,z)u_2(a,z),\\
    V(a,z)&=V(a,0)u_1(a,z)+V'(a,z)u_2(a,z),
\end{align}
where 
\begin{align}
    u_1(a,z)&=e^{-\frac{1}{4}z^2}\left[1+\sum_{p=1}^{\infty}\left(\prod_{k=1}^{p}(a-\frac{3}{2}+2k)\frac{z^{2p}}{(2p)!}\right)\right],\\
    u_2(a,z)&=e^{-\frac{1}{4}z^2}\left[z+\sum_{p=1}^{\infty}\left(\prod_{k=1}^{p}(a-\frac{1}{2}+2k)\frac{z^{2p+1}}{(2p+1)!}\right)\right],
\end{align} 
and
\begin{align}
    U(a,0) &=\frac{\sqrt{\pi}}{2^{\frac{a}{2}+\frac{1}{4}}\Gamma(\frac{3}{4}+\frac{a}{2})}, \label{eq:U0} \\
    U'(a,0)&=-\frac{\sqrt{\pi}}{2^{\frac{a}{2}-\frac{1}{4}}\Gamma(\frac{1}{4}+\frac{a}{2})}, \label{eq:Up0} \\
    V(a,0) &=\frac{\pi 2^{\frac{a}{2}+\frac{1}{4}}}{\left(\Gamma(\frac{3}{4}-\frac{a}{2})\right)^2\Gamma(\frac{1}{4}+\frac{a}{2})}, \label{eq:V0} \\
    V'(a,0)&=\frac{\pi 2^{\frac{a}{2}+\frac{3}{4}}}{\left(\Gamma(\frac{1}{4}-\frac{a}{2})\right)^2\Gamma(\frac{3}{4}+\frac{a}{2})}. \label{eq:Vp0}
\end{align}
For $z$ large, $a$ moderate, and $z \gg |a|$, the asymptotic expansions of $U(a,z)$ and $V(a,z)$ are given by~\cite{zhang1996computation}:
\begin{align}
    U(a,z)&\sim e^{-\frac{1}{4}z^{2}}z^{-a-\frac{1}{2}}\left[1-\frac{(a+\frac{1}{2})(a+\frac{3}{2})}{2z^2}+\ldots\right],\\
    V(a,z)&\sim \sqrt{\frac{2}{\pi}} \,e^{\frac{1}{4}z^{2}}z^{a-\frac{1}{2}}\left[1+\frac{(a-\frac{1}{2})(a-\frac{3}{2})}{2z^2}+\ldots\right].
    \label{eq:asymptotic}
\end{align}
Using Eq.~\eqref{eq:Uarg} we have:
\begin{align*}
U(a,z)+U(a,-z)&=U(a,z)-\sin(\pi a) U(a,z)+\frac{\pi}{\Gamma(\frac{1}{2}+a)}V(a,z)\\
              &=\left[1-\sin(\pi a)\right]U(a,z)+\frac{\pi}{\Gamma(\frac{1}{2}+a)}V(a,z).
\end{align*}
Using Eq.~\eqref{eq:DUrel} and the above equation, we get: 
\begin{eqnarray}
D_{\beta-1}(-\sqrt{2}Z_i)+D_{\beta-1}(\sqrt{2}Z_i) & = & U\left(\frac12-\beta,-\sqrt{2}Z_i\right)+U\left(\frac12-\beta,\sqrt{2}Z_i\right) \nonumber \\
& = & g(\beta)U\left(\frac12-\beta,\sqrt{2}Z_i\right)+ h(\beta)V\left(\frac12-\beta,\sqrt{2}Z_i\right),
\label{eq:uvform}
\end{eqnarray}
where $g(\beta)=\left[1-\sin((\frac12-\beta)\pi)\right]$ and $h(\beta)=\frac{\pi}{\Gamma(1-\beta)}$.

\begin{figure}[htb]
   \centering\includegraphics[width=0.5\textwidth]{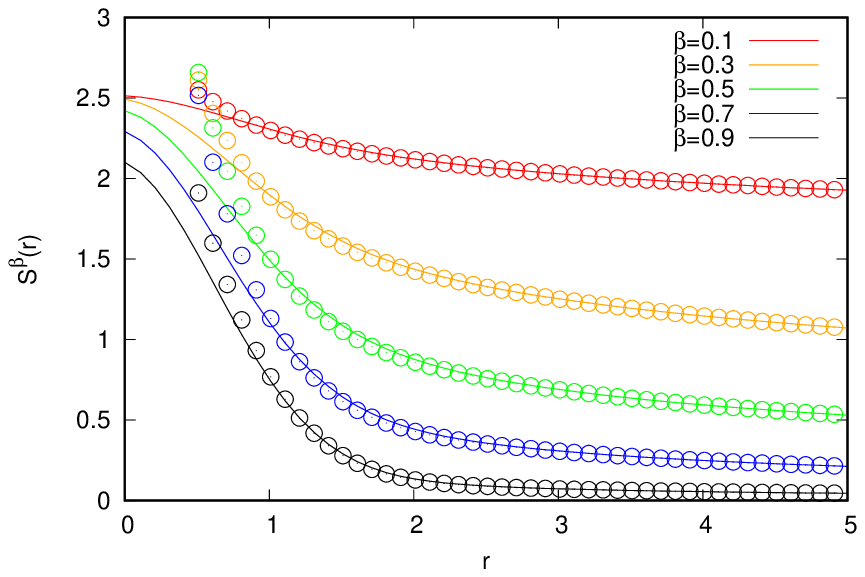}
    \caption{$S^\beta(r)$ from Eq.~\eqref{eq:Sdef}; symbols are referred to the asymptotic approximation.
    \label{fig:plots}}
\end{figure}

The parabolic cylinder function $D_\nu(z)$ can be readily evaluated using the routine \verb+pbdv+, available within the FORTRAN77 code from Ref.~\cite{zhang1996computation}.
A FORTRAN90 version is available~\cite{specialfuncs}, and has been used to produce the illustrations of Fig.~\ref{fig:plots}.
For large arguments, we can rely on the asymptotic expansions of $U$ and $V$, namely approximate terms of the form
\begin{equation}
S^\nu(z) \equiv e^{\frac{-z^2}{2}}\left(D_{\nu-1}(-\sqrt{2}z)+D_{\nu-1}(\sqrt{2}z)\right), \label{eq:Sdef}
\end{equation}
which appear in the particle representation of the function $\tilde u(x,t)$, see Eq.~\eqref{eq:utilde_def} and~\eqref{eq:utilde_exp},
and also in estimates of the fractional diffusion flux and of the fractional diffusion term.
Specifically, for $\nu=\beta$, $0 < \beta < 1$, we have
\begin{equation}\label{eq:g_asymp}
   S^\beta (z) = e^{\frac{-z^2}{2}}\left(g(\beta)U\left(\frac12-\beta,\sqrt{2}z\right)+ h(\beta)V\left(\frac12-\beta,\sqrt{2}z\right)\right), 
\end{equation}
where from Eq.~\eqref{eq:asymptotic}:
\begin{equation}
\begin{aligned}
U(\frac12-\beta,z)\sim e^{-\frac{1}{4}z^{2}}z^{\beta-1}\left[1-\frac{(1-\beta)(2-\beta)}{2z^2}\right]\\
    V(\frac12-\beta,z)\sim \sqrt{\frac{2}{\pi}} \,e^{\frac{1}{4}z^{2}}z^{-\beta}\left[1+\frac{\beta(\beta+1)}{2z^2}\right].
\end{aligned}
\end{equation}
As a consequence, $S^\beta(z)\propto z^{-\beta}$ for $z\to\infty$.
Figure~\ref{fig:plots} shows $S^\beta(x)$, which is represented both in exact form and in asymptotic form truncated at the first order.


\end{document}